\documentclass[envcountsame,envcountsect,runningheads]{cl2emult}

\usepackage{amsmath,amssymb,amsfonts,latexsym,epsfig}
\usepackage[english,frenchb]{babel}
\usepackage[T1]{fontenc} 
\usepackage{enumerate}
\usepackage{eucal}

\usepackage{makeidx}  
\usepackage{graphicx} 
\usepackage{subeqnar} 
\usepackage{multicol} 
\usepackage{cropmark} 
\usepackage{amscd}
\makeindex            
\begin{document}

\frontmatter

\thispagestyle{empty}
{\LARGE\parindent=0pt
{\bfseries Marc Burger\\Alessandra Iozzi\\\par\large Editors
\par\vspace{2cm}
\rule{\textwidth}{1pt}
{\itshape \LARGE RIGIDITY IN DYNAMICS\\ AND GEOMETRY\\
\par\large Contributions from the Programme ``Ergodic Theory,\\
Geometric Rigidity and Number Theory''\\
\par Isaac Newton Institute for the Mathematical Sciences\\
Cambridge, United Kingdom\\
\par 5 January -- 7 July 2000}\\
\rule{\textwidth}{1pt}}}
\vfill
{\large\itshape
\begin{tabular}{@{}l}
{\Large\rm\bfseries Springer}\\[8pt]
Berlin\enspace Heidelberg\enspace New\kern0.1em York\\[5pt]
Barcelona\enspace Budapest\enspace Hong\kern0.2em Kong\\[5pt]
London\enspace Milan\enspace Paris\enspace\\[5pt]
Santa\kern0.2em Clara\enspace Singapore\enspace Tokyo
\end{tabular}}
\clearpage

\selectlanguage{english}
\setcounter{tocdepth}{1}

\mainmatter
%
\spnewtheorem*{literature}{Pointers to the
literature}{\itshape}{\rmfamily}


  \newcounter{CASE}
  \newenvironment{CASE}[1][\unskip]{\refstepcounter{CASE}\em
  \medbreak
   \noindent Case \theCASE\ #1.\ }{\unskip\upshape}
  \renewcommand{\theCASE}{\arabic{CASE}}

  \newcounter{STEP}
\newenvironment{STEP}[1][\unskip]{\refstepcounter{STEP}
\em
   \medbreak\noindent Step \theSTEP\
  #1.\ }{\unskip\upshape}


\newcounter{ZarEgCount}

\newcommand{\ZarEg}{\stepcounter{ZarEgCount}
  \arabic{ZarEgCount}) }


\newcommand{\pref}[1]{{\upshape(}\ref{#1}{\upshape)}}
\newcommand{\mysee}[1]{{\upshape(}see~\ref{#1}{\upshape)}}
\newcommand{\fullref}[2]{\ref{#1}\pref{#1-#2}}
\newcommand{\fullsee}[2]{{\upshape(}see~\fullref{#1}{#2}{\upshape)}}


\renewcommand{\operatorname}[1]{\mathop{\rm#1}}
\newcommand{\graph}{\operatorname{graph}}
\newcommand{\closure}{\overline}
\newcommand{\cover}{\widetilde}
\newcommand{\torus}{\mathbb{T}}
\newcommand{\Aut}{\operatorname{Aut}}

\newcommand{\SO}{\operatorname{SO}}
\newcommand{\SL}{\operatorname{SL}}
\newcommand{\GL}{\operatorname{GL}}
\newcommand{\SU}{\operatorname{SU}}
\newcommand{\Id}{\operatorname{Id}}
\newcommand{\real}{\mathord{\mathbb{R}}}
\newcommand{\complex}{\mathord{\mathbb{C}}}
\newcommand{\integer}{\mathord{\mathbb{Z}}}
\newcommand{\Rrank}{\operatorname{\hbox{\upshape$\real$-rank}}}
\newcommand{\iso}{\cong}
\newcommand{\diffeo}{\simeq}
\newcommand{\Ad}{\operatorname{Ad}\nolimits}

\newcommand{\bigset}[2]{\left\{\, #1
  \mathrel{\left| \vphantom {\left\{ #1 \mid #2 \right\} } \right.}
  #2 \,\right\} }

\title*{Superrigid Subgroups
  and Syndetic Hulls\protect\newline
  in Solvable Lie Groups}

\toctitle{Superrigid Subgroups and Syndetic Hulls\protect\newline
  in Solvable Lie Groups}

\titlerunning{Superrigid Subgroups and Syndetic Hulls
  in Solvable Lie Groups}

\author{Dave Witte\thanks{The preparation of this paper was partially supported by a
grant from the National Science Foundation (DMS-9801136).}}

\institute{Department of Mathematics,
  Oklahoma State University,
  Stillwater, OK 74078, USA\\
 {\em e-mail:} {\tt dwitte@math.okstate.edu}}

\setcounter{page}{441}   

\maketitle

\begin{abstract}
  It is not difficult to see that every group homomorphism
from~$\integer^k$ to~$\real^n$ extends to a homomorphism
from~$\real^k$ to~$\real^n$. We discuss other examples of
discrete subgroups~$\Gamma$ of connected Lie groups~$G$, such
that the homomorphisms defined on~$\Gamma$ can (``virtually'')
be extended to homomorphisms defined on all of~$G$. For the
case where $G$ is solvable, we give a simple proof that
$\Gamma$ has this property if it is Zariski dense. The key
ingredient is a result on the existence of syndetic hulls.
  \end{abstract}

\section{What Is a Superrigid Subgroup?}

Let us begin with a trivial example of the type of theorem
that we will discuss. It follows easily from the fact
that a linear transformation can be defined to have any
desired action on a basis. (See Sect.~\ref{superproofSect}
for a more complicated proof.)

\begin{proposition} \label{ZktoRd}
  Any group homomorphism $\varphi  \colon \integer^k \to \real ^d$
extends to a continuous homomorphism $\hat\varphi  \colon \real^k
\to \real ^d$.
  \end{proposition}

A superrigidity theorem is a version of this simple proposition
in the situation where $\integer^k$, $\real^k$, and~$\real^d$
are replaced by more interesting groups. Suppose
  $\Gamma$ is a discrete subgroup of a connected Lie group~$G$,
and
  $H$ is some other Lie group.
  Does every homomorphism $\varphi \colon \Gamma \to H$ extend to
a continuous homomorphism $\hat\varphi$ defined on all of~$G$?

All of the Lie groups we consider are assumed to be
\emph{linear groups}\index{group!linear}; that is, they are
subgroups of $\GL(\ell,\complex)$, for some~$\ell$. For example,
$\real^d$ can be thought of as a linear group; in particular:
  \begin{equation} \label{R3inGL4}
  \real^3 \iso
  \begin{pmatrix}
  1&0&0&\real\cr
  0&1&0&\real\cr
  0&0&1&\real\cr
  0&0&0&1 \cr
  \end{pmatrix}
   \; .
  \end{equation}
  Thus, any homomorphism into $\real^d$ can be thought of as a
homomorphism into $\GL(d+1,\real)$. The study of homomorphisms
into $\GL(d,\real)$ or $\GL(d,\complex)$ is known as
\emph{Representation Theory}. Unfortunately, in this much more
interesting setting, not all homomorphisms extend.

\begin{proposition}
  There is a group homomorphism $\varphi \colon \integer \to
\GL(d,\real)$ that does not have a continuous extension to a
homomorphism $\hat\varphi  \colon \real \to \GL(d,\real)$.
  \end{proposition}

\begin{proof}
  Fix a matrix $A \in \GL(d,\real)$, such that
  $\det A < 0$,
  and define $\varphi(n) = A^n$, for $n \in
\integer$. Then $\varphi(m+n) = \varphi(m) \cdot \varphi(n)$, so $\varphi$
is a group homomorphism. Since $\det X \neq 0$, for all $X \in
\GL(d,\real)$, there is no continuous function $\hat\varphi \colon
\real \to \GL(d,\real)$, such that $\hat\varphi(0) = \Id$ and
$\hat\varphi(1) = A$. Hence, $\varphi$ does not have a continuous
extension to~$\real$.
  \qed \end{proof}

The example shows that we cannot expect $\hat\varphi(n)$ to equal
$\varphi(n)$ for \emph{all} $n \in \integer$, so we relax the
restriction to require equality only when $n$ belongs to some
finite-index subgroup of~$\integer$. In group theory, it is
standard practice to say that a group \emph{virtually} has a
property if some finite-index subgroup has the property. In
that spirit, we make the following definition.

\begin{definition}
  Suppose $\Gamma$ is a subgroup
of~$G$. We say that a homomorphism $\varphi \colon \Gamma \to H$
\emph{virtually extends}\index{virtually!extends} to a
homomorphism $\hat \varphi \colon G \to H$ if there is a
finite-index subgroup~$\Gamma'$ of~$\Gamma$, such that
$\varphi(\gamma) = \hat\varphi(\gamma)$ for all $\gamma \in \Gamma'$.
  \end{definition}

The following result is not as trivial as
Proposition~\ref{ZktoRd}.

\begin{proposition} \label{ZktoGL}
  Any group homomorphism $\varphi \colon \integer^k \to
\GL(d,\real)$ virtually extends to a continuous homomorphism
$\hat\varphi  \colon \real^k \to \GL(d,\real)$. Similarly for
homomorphisms into $\GL(d,\complex)$.
  \end{proposition}

Proposition~\ref{ZktoGL} has the serious weakness that it
gives no information at all about the image of~$\hat\varphi$. A
superrigidity theorem should state not only that a virtual
extension exists, but also that if the image of the original
homomorphism~$\varphi$ is well behaved, then the image of the
extension~$\hat\varphi$ is similarly well behaved. For example,
we want:
  \begin{itemize}
  \item If $\hat \varphi(\Gamma) \subset \real^d$, then
$\hat\varphi(G) \subset \real^d$ (as embedded in~\pref{R3inGL4}).
  \item If all the matrices in $\varphi(\Gamma)$ commute with each
other, then all the matrices in
$\hat\varphi(G)$ commute with each other.
  \item If all of the matrices in $\varphi(\Gamma)$ are upper
triangular, then all of the matrices
in $\hat\varphi(G)$ are upper triangular.
  \end{itemize}
  All of these properties, and many more, are obtained by
requiring that $\hat\varphi(G)$ be contained in the ``Zariski
closure'' of~$\varphi(\Gamma)$.

The Zariski closure will be formally defined in
Sect.~\ref{SolvZarSect}. For now, it suffices to have
an intuitive understanding:
  \begin{quote}
  The Zariski closure\index{Zariski!closure}~$\closure{\Gamma}$ of
a subgroup~$\Gamma$ of $\GL(\ell,\complex)$ is the ``natural''
virtually connected subgroup of $\GL(\ell,\complex)$ that
contains~$\Gamma$.
  \end{quote}
  (By ``virtually connected\index{connected!virtually}\index{virtually!connected},'' 
we mean that the Zariski closure, although perhaps not connected,
has only finitely many components.) Some examples should help to
clarify the idea.

\begin{example} \label{ZarEgs}
  \setcounter{ZarEgCount}{0}
  \ZarEg $\real^d$, as embedded in~\pref{R3inGL4} is its own
Zariski closure; we say that $\real^d$ is \emph{Zariski
closed}\index{Zariski!closed}.

\ZarEg Let
  $$
  G_1 =
  \begin{pmatrix}
  1&\real&\real&\real\cr
  0&1&\real&\real\cr
  0&0&1&\real\cr
  0&0&0&1\cr
  \end{pmatrix}
  \mbox{\qquad and\qquad }
  \Gamma_1  = \begin{pmatrix}
  1&\integer&\integer&\integer\cr
  0&1&\integer&\integer\cr
  0&0&1&\integer\cr
  0&0&0&1\cr
  \end{pmatrix}
  \; .
  $$
  Then $G_1$ is a perfectly natural, connected subgroup, so
$G_1$ is Zariski closed. Because $G_1$ is the natural
connected subgroup that contains~$\Gamma_1$, we have
  $ \closure{\Gamma_1} = G_1 $.
  (We may say that $\Gamma_1$ is \emph{Zariski
dense}\index{Zariski!dense} in~$G_1$.)

\ZarEg Let
  $$G_2 =
  \begin{pmatrix}
  1&\real&\complex\cr
  0&1&0\cr
  0&0&1\cr
  \end{pmatrix}
  \mbox{\qquad and\qquad }
  \Gamma_2 =
  \begin{pmatrix}
  1&\integer&\integer +\integer i\cr
  0&1&0 \cr
  0&0&1\cr
  \end{pmatrix}
  \; .$$
  Then
  $\closure{\Gamma_2} = G_2 = \closure{G_2}$.

\ZarEg Let
  \begin{equation} \label{G3Defn}
  G_3 =
  \bigset{
  \begin{pmatrix}
  1&t&z\cr
  0&1&0\cr
  0&0&e^{2\pi it}\cr
  \end{pmatrix}
  }{
  \begin{matrix}
  t \in \real, \cr
  z \in \complex
  \end{matrix}
  }
  \; .
  \end{equation}
  Then, although $G_3$ is connected, it is not Zariski closed. The
notion of Zariski closure comes from Algebraic Geometry, where
only polynomial functions are considered. Thus, because the
exponential function is transcendental, not polynomial, an
Algebraic Geometer does not see the coupling between the $(1,2)$
entry and the $(3,3)$ entry of the matrix; so, from an Algebraic
Geometer's point of view, there is no constraint linking these
two matrix entries. The $(1,2)$ entry takes any real value, the
$(3,3)$ entry takes any value on the unit circle, and the Zariski
closure allows these values entirely independently:
  $$ \closure{G_3} =
   \begin{pmatrix}
  1&\real&\complex\cr
  0&1&0\cr
  0&0&\torus\cr
  \end{pmatrix}
  \; .$$
  (The analogous example for a topologist would be a discontinuous
function $f \colon \real \to \torus$, such that the graph of~$f$
is dense in $\real \times \torus$.)
  As another important observation, note that $\Gamma_2 \subset
G_3$. However, we know that $\closure{\Gamma_2} = G_2  \neq
\closure{G_3}$, so $\Gamma_2$ is \emph{not} Zariski dense
in~$G_3$.
  \end{example}

We can now define a version of superrigidity:

\begin{definition}
  A discrete subgroup~$\Gamma$ of a Lie group~$G$ is
\emph{superrigid}\index{subgroup!superrigid} in~$G$ if every
homomorphism $\varphi \colon \Gamma \to \GL(d,\real)$ virtually
extends to a continuous homomorphism $\hat\varphi \colon G \to
\closure{\varphi(\Gamma)}$.
  \end{definition}

The following superrigidity result (a special case of the main
theorem stated later in this section) strengthens
Proposition~\ref{ZktoGL}. Except for the minor discrepancy
between extensions and virtual extensions, it also generalizes
Proposition~\ref{ZktoRd}.

\begin{proposition} \label{Zksuper}
  $\integer^k$ is superrigid in~$\real^k$.
  \end{proposition}

A more interesting result (one that deserves to be
called a theorem) applies to nonabelian groups. In this
section, we consider only solvable groups:

\begin{definition} \label{solv<>triang}
  A connected Lie subgroup~$G$ of $\GL(\ell,\complex)$ is
\emph{solvable}\index{group!solvable Lie}\index{Lie!group!solvable} if (perhaps after a
suitable change of basis) it is upper triangular.

(This is not the usual definition, but it is more concrete, and
it is equivalent to the usual one in our setting
\mysee{SolvDefnsSame}.)
  \end{definition}

  For example, the groups~$G_1$, $G_2$, and~$G_3$ defined in
Example~\ref{ZarEgs} are obviously solvable. Also, note that
any set of pairwise commuting matrices can be simultaneously
triangularized, so abelian groups are solvable.

To avoid technical problems that could force us to pass to a
finite cover, we usually assume that the fundamental group of~$G$
is trivial:

\begin{definition}
  A Lie group~$G$ is \emph{$1$-connected}\index{connected!$1$---} if it
is connected and simply connected.
  \end{definition}

The following example shows that, for some solvable groups, not
all subgroups are superrigid.

\begin{example}  \label{notsuper}
  Let
  $$
  G =
  \begin{pmatrix}
  \real^+ & \real \cr
  0 & 1
  \end{pmatrix}
  \mbox{\qquad and\qquad }
  \Gamma =
  \begin{pmatrix}
  1 & \integer \cr
  0 & 1
  \end{pmatrix}
  \; .$$
  Then $G$ is obviously solvable, and $\Gamma$ is a discrete
subgroup.

Any homomorphism $G \to \real$ must vanish on the commutator
subgroup
  $$[G,G] =
  \begin{pmatrix}
  1 & \real \cr
  0 & 1
  \end{pmatrix}
  \supset \Gamma
  \; ,$$
  so the only homomorphism $\varphi \colon \Gamma \to \real$
that virtually extends to~$G$ is the trivial homomorphism.
(Because $\real$ has no nontrivial finite subgroups, any
virtually trivial homomorphism into~$\real$ must actually be
trivial.) Therefore, not all homomorphisms virtually extend, so
$\Gamma$ is not superrigid in~$G$.
  \end{example}

The moral of this example is that small subgroups cannot be
expected to be superrigid: if $\Gamma$ is only a small part
of~$G$, then a homomorphism defined on~$\Gamma$ knows nothing
about most of~$G$, so it cannot be expected to be compatible
with the structure of all of~$G$. This suggests that, to
obtain a superrigidity theorem, we should assume that $\Gamma$
is large, in some sense. The correct sense is Zariski density.

\begin{theorem}[Witte 1997] \label{solvsuper}
  If $\Gamma$ is a discrete subgroup of a $1$-connected,
solvable Lie group $G \subset \GL(\ell,\complex)$, such that
$\closure{\Gamma} = \closure{G}$, then $\Gamma$ is
superrigid\index{subgroup!superrigid} in~$G$.
  \end{theorem}

\begin{example}
  Because
$\closure{\Gamma_1} = \closure{G_1}$ and $\closure{\Gamma_2} =
\closure{G_2}$, the theorem implies that $\Gamma_1$ is
superrigid in~$G_1$, and $\Gamma_2$ is superrigid in $G_2$.
(Actually, the case of~$G_2$ follows already from
Proposition~\ref{Zksuper}, because it is easy to see that $G_2
\iso \real^3$.) Because $\integer^k$ is Zariski dense
in~$\real^k$ (for $\real^k$ as in \pref{R3inGL4}),
Proposition~\ref{Zksuper} is a special case of this theorem.

On the other hand, we have $\closure{\Gamma_2} \neq
\closure{G_3}$, and, although the theorem does not tell us
this, it is easy to see that $\Gamma_2$ is not superrigid
in~$G_3$. (The subgroup~$\Gamma_2$ is abelian, and the
intersection $\Gamma_2 \cap [G_3,G_3]$ is infinite, so this is
much the same as Example~\ref{notsuper}.)
  \end{example}

Unfortunately, Proposition~\ref{ZktoRd} is not quite a
corollary of this theorem, because of the discrepancy between
a virtual extension and an actual extension.
Section~\ref{othersuperSect} states a version of
Theorem~\ref{solvsuper} that, under additional technical
hypotheses, provides an actual extension, thus generalizing
Proposition~\ref{ZktoRd}. The section also states a more
precise version of Theorem~\ref{solvsuper} that determines
exactly which subgroups of a solvable Lie group are superrigid,
and briefly discusses superrigidity theorems for Lie groups
that are not solvable.

A simple proof of Theorem~\ref{solvsuper} will be given in
Sect.~\ref{superproofSect}, modulo an assumption about the
existence of syndetic hulls. This gap will be filled in
Sect.~\ref{SyndeticSect}, after some definitions and basic
results are recalled from the literature in
Sect.~\ref{SolvZarSect}.

\section{Other Superrigidity Theorems} \label{othersuperSect}

Theorem~\ref{solvsuper} can be extended to stronger results that
provide more detailed information about solvable groups, and to
broader results that apply to more general groups.

\subsection{More on Superrigid Subgroups of Solvable Groups}

  There are two obvious reasons that the converse of
Theorem~\ref{solvsuper} does not hold.
  \begin{itemize}
  \item If $\Gamma$ is superrigid in~$B$, then $e \times
\Gamma$ is superrigid in $A \times B$. So, to be superrigid in
a direct product, it suffices to be Zariski dense in one of
the factors. This generalizes to semidirect products, as well.
  \item The group~$G$ has many different representations in
$\GL(\ell,\complex)$; it may happen that $\Gamma$ is Zariski
dense in some of these embeddings, but not in others. (For
example, if we realize $\real$ as the subgroup of~$G_3$ with
$z = 0$, then $\integer$ is not Zariski dense in~$\real$.) This
ambiguity is eliminated by using only the adjoint
representation (even though this is not an embedding if $G$
has a center).
  \end{itemize}
  The following corollary shows that these two obvious reasons
are the only ones.

\begin{corollary} \label{solvsuperiff}
  A discrete subgroup~$\Gamma$ of a $1$-connected,
solvable Lie group~$G$ is superrigid if and
only if
  \begin{enumerate}
  \item $G = A \rtimes B$, for some closed subgroups~$A$
and~$B$ of~$G$; such that
  \item $B$ contains a finite-index subgroup $\Gamma'$
of~$\Gamma$; and
  \item $\closure{\Ad_B\Gamma'} = \closure{\Ad_B B}$.
  \end{enumerate}
  \end{corollary}

\begin{proof}[$\Rightarrow$] Let $B = \closure{\Gamma}^\circ$. The
inclusion $\varphi \colon \Gamma \hookrightarrow B$ (virtually)
extends to a continuous homomorphism $\hat\varphi \colon G \to
B$. Since $\hat\varphi|_\Gamma = \Id_\Gamma$, and $\Gamma$ is
Zariski dense in~$B$, it is reasonable to expect that
$\hat\varphi|_B = \Id_B$. (Actually, this need not quite be true,
but it is close to correct.) Then $G = (\ker\hat\varphi) \rtimes
B$.
  \qed \end{proof}

The most important special case is when $\Gamma$ is a lattice
in~$G$:

\begin{definition}
  A discrete subgroup~$\Gamma$ of a solvable Lie group~$G$ is a
\emph{lattice}\index{lattice} if $G/\Gamma$ is compact.
  \end{definition}

For example, $\integer^k$ is a lattice in~$\real^k$, and, in
Example~\ref{ZarEgs}, $\Gamma_1$ is a lattice in~$G_1$, and
$\Gamma_2$ is a lattice in both~$G_2$ and~$G_3$. The
superrigidity criterion for lattices is very simple:

\begin{corollary} \label{solvlattsuper}
  A lattice~$\Gamma$ in a $1$-connected, solvable Lie
group~$G$ is superrigid if and only if
  $\closure{\Ad_G\Gamma} = \closure{\Ad_G G}$.
  \end{corollary}

The following result provides an extension, not just a virtual
extension, under mild hypotheses on~$\varphi$.

\begin{corollary} \label{noneedfinind}
  Let $\Gamma$ be a lattice in a connected, solvable Lie
group~$G$, such that $\closure{\Ad_G\Gamma} = \closure{\Ad_G
G}$. If $\varphi \colon \Gamma \to \GL(d,\complex)$ is a
homomorphism, such that
  \begin{itemize}
  \item  $\varphi(\Gamma) \subset \closure{\varphi(\Gamma)}^\circ$,
  \item the center of $\closure{\varphi(\Gamma)}^\circ$ is
connected, and
  \item $\varphi \bigl( \Gamma \cap [G,G] \bigr)$ is unipotent,
  \end{itemize}
  then $\varphi$ extends to a continuous homomorphism $\hat\varphi
\colon G \to \closure{\varphi(\Gamma)}$.
  \end{corollary}

The groups $G_2$ and~$G_3$ of Example~\ref{ZarEgs} are
non-isomorphic solvable groups that have isomorphic lattices
(namely,~$\Gamma_2$). The following consequence of
superrigidity implies that solvable groups with isomorphic
lattices differ only by rotations being added to and/or
removed from their Zariski closures.

\begin{corollary} \label{solvlattiso}
  Suppose
  \begin{itemize}
  \item $\Gamma _1$ and~$\Gamma _2$ are lattices in
$1$-connected, solvable Lie groups $G_1$ and~$G_2$,
  \item $\closure{\Ad_{G_1} \Gamma_1} = \closure{\Ad_{G_1}
G_1}$, and
  \item $\pi \colon \Gamma _1 \to \Gamma _2$ is an isomorphism.
  \end{itemize}
  Then $\pi$ extends to an embedding $\sigma \colon G_1 \to T
\ltimes G_2$, for any maximal compact subgroup~$T$
of~$\closure{\Ad_{G_2} G_2}^\circ$.
  \end{corollary}

\begin{corollary}[Mostow 1954] \label{SolvMostow}
  Suppose
  \begin{itemize}
  \item $\Gamma _1$ and~$\Gamma _2$ are lattices in
$1$-connected, solvable Lie groups $G_1$ and~$G_2$, and
  \item $\Gamma_1$ is isomorphic to $\Gamma_2$.
  \end{itemize}
  Then $G_1/\Gamma _1$ is diffeomorphic to $G_2/\Gamma _2$.
  \end{corollary}

\subsection{Superrigid Subgroups of Semisimple Groups}

For groups that are not solvable, both ``lattice'' and
``superrigid'' need to be generalized from the definitions
above.

\begin{definition}
  \begin{itemize}
  \item A discrete subgroup~$\Gamma$ of a Lie group~$G$ is a
\emph{lattice}\index{lattice} if there is $G$-invariant Borel
probability measure on $G/\Gamma$.

\item A lattice $\Gamma$ in a Lie group~$G$ is
\emph{superrigid}\index{subgroup!superrigid} if, for every
homomorphism $\varphi \colon \Gamma \to \GL(d,\complex)$, there is a
compact, normal subgroup~$K$ of $\closure{\varphi(\Gamma)}^\circ$, a
continuous homomorphism $\hat\varphi \colon G \to
\closure{\varphi(\Gamma)}^\circ/K$, and a finite-index
subgroup~$\Gamma'$ of~$\Gamma$, such that
  $\hat\varphi(\gamma) = \varphi(\gamma) K$, for all $\gamma \in
\Gamma'$.
  \end{itemize}
  \end{definition}

The semisimple case is orders of magnitude more difficult than
the solvable case. We still do not have a complete answer, but
the following amazing theorem of G.~A.~Margulis settles most
cases.

\begin{theorem}[Margulis Superrigidity Theorem]
\label{MargSuper}
  If $n \ge 3$, then every lattice in $\SL(n,\real)$ is
superrigid.\index{theorem!Margulis!superrigidity}\index{superrigidity!Margulis --- theorem}\index{Margulis!theorem!superrigidity}

The same is true for irreducible lattices in any other
connected, semisimple, linear Lie group~$G$ with $\Rrank G \ge
2$.
  \end{theorem}

K.~Corlette \cite{Corlette-super} proved that lattices in
$\operatorname{Sp}(1,n)$ are superrigid, and also lattices in
the exceptional group of real rank one. Thus, to complete the
study of lattices in semisimple groups, all that remains is to
determine which lattices in $\SO(1,n)$ and $\SU(1,n)$ are superrigid.
(Many lattices in $\SO(1,n)$ are \emph{not} superrigid.)
Lattices are not the whole story, however: H.~Bass and
A.~Lubotzky \cite{BassLubotzky} recently constructed an
example of a Zariski dense superrigid discrete subgroup~$\Gamma$ of a
semisimple group, such that $\Gamma$ is \emph{not} a lattice.

A superrigidity theorem describes a very close connection
between a lattice~$\Gamma$ and the ambient Lie group~$G$. In
fact, for semisimple groups, the connection is so close that
superrigidity tells us almost exactly what the lattice must be.
In all of our examples above, the lattice~$\Gamma$ consists of
the integer points of~$G$. The following major consequence of
the Margulis Superrigidity Theorem implies that this is
essentially the only way to make a lattice in a simple group
of higher real rank. (However, one needs to allow certain
algebraic integers in place of ordinary integers.)

\begin{corollary}[Margulis Arithmeticity Theorem]
\label{MargArith}
  If $n \ge 3$, then every lattice in $\SL_n(\real)$ is
arithmetic.\index{theorem!Margulis!arithmeticity}\index{arithmeticity!Margulis --- theorem}\index{Margulis!theorem!arithmeticity}

The same is true for irreducible lattices in any connected,
semisimple, linear Lie group~$G$ with $\Rrank G \ge 2$.
  \end{corollary}

\subsection{Superrigid Subgroups of Other Lie Groups}

The following proposition makes it easy to combine the
semisimple case with the solvable case.

\begin{proposition}[L.~Auslander] \label{Auslander-compatible}
  Let $G = R \rtimes L$ be a Levi decomposition of a connected
Lie group~$G$, and let $\sigma \colon G \to L$ be the
corresponding quotient map. If $\Gamma$ is a lattice in~$G$,
such that $\closure{\Ad_G \Gamma} = \closure{\Ad_G G}$, then
$\Gamma \cap R$ is a lattice in~$R$, and $\sigma(\Gamma)$ is a
lattice in~$L$.
  \end{proposition}

\begin{corollary} \label{supermixed}
  Let $G = R \rtimes L$ be a Levi decomposition of a connected,
linear Lie group~$G$, and let $\sigma \colon G \to L$ be the
corresponding quotient map.
  A lattice $\Gamma$ in~$G$ is superrigid if and only if
   \begin{itemize}
  \item there is a compact, normal subgroup~$C$ of
$\closure{\Ad_G G}$, such that $\bigl( \closure{\Ad_G \Gamma}
\bigr) C =\closure{\Ad_G G}$, and
  \item the lattice $\sigma(\Gamma)$ is superrigid in~$L$.
  \end{itemize}
  \end{corollary}

\begin{literature}
  Corollary~\ref{solvsuperiff} is from \cite{Witte-SolvDense}.
Corollaries~\ref{solvlattsuper}, \ref{noneedfinind},
\ref{solvlattiso}, and~\ref{supermixed} are from
\cite{Witte-SolvRig}. (See \cite{Starkovwi} for results related
to~\ref{solvlattiso}, without the compact subgroup~$T$.)
Corollary~\ref{SolvMostow} and
Proposition~\ref{Auslander-compatible} appear in
\cite[Theorems~3.6 and~8.24]{Raghunathan?}.
Theorems~\ref{MargSuper} and Corollary~\ref{MargArith} are
discussed in \cite{MargulisBook} and~\cite{ZimmerBook}.
  \end{literature}

\section{Our Prototypical Proof of Superrigidity}
\label{superproofSect}

We now give a proof of Proposition~\ref{ZktoRd} that is
somewhat more difficult than necessary, because this argument
can be generalized to other groups.

\begin{proof}[of Proposition~\ref{ZktoRd}]
  Let
  \begin{itemize}
  \item
  $\hat\Gamma = \graph(\varphi)
  = \bigl\{\, \bigl( \gamma, \varphi(\gamma) \bigr)
  \mid \gamma \in \Gamma \, \bigr\}
  \subset \real^k \times \real^d $,
  \item $ X = \operatorname{span} \hat\Gamma $
  be the subspace of the vector space~$\real^k \times \real^d$
spanned by~$\hat\Gamma$, and
  \item $p \colon \real^k \times \real^d \to \real^k$ be the
natural projection onto the first factor.
  \end{itemize}

\setcounter{STEP}{0}

\begin{STEP} \label{superabelPf-ponto}
  We have $p(X) = \real^k$.
  \end{STEP}
  Note that:
  \begin{itemize}
  \item $p(X)$ is connected (because $X$ is connected and $p$ is
continuous);
  \item $p(X)$ is an additive subgroup of~$\real^k$ (because
$X$ is an additive subgroup, and $p$ is an additive
homomorphism); and
  \item $p(X)$ contains $\integer^k$ (because $p(X)$ contains
$p(\hat\Gamma) = \operatorname{dom}(\varphi) = \integer^k$).
  \end{itemize}
  Since
  \begin{equation} \label{abelian-fullsubgroup}
  \mbox{no connected, proper subgroup of~$\real^k$
contains~$\integer^k$,}
  \end{equation}
  the desired conclusion follows.

\begin{STEP} \label{superabelPf-kerp}
  We have $X \cap (0 \times \real^d) = 0$.
  \end{STEP}
  Because $\Gamma$ is discrete, we know that $\varphi$ is
continuous, so the $\varphi$-image of any compact subset
of~$\Gamma$ is compact. This implies that $p|_{\hat\Gamma}$,
the restriction of~$p$ to~$\hat\Gamma$, is a proper map. (That
is, the inverse image of every compact set is compact.)
  It is a fact that
  \begin{equation} \label{abelian-syndetic}
  \mbox{$(\operatorname{span} \Lambda)/\Lambda\;$ is compact, for
every closed subgroup $\Lambda$ of $\real^k \times \real^d$;}
  \end{equation}
  therefore, $X = \operatorname{span} \hat\Gamma$ differs
from~$\hat\Gamma$ by only a compact amount.
  Since $p|_{\hat\Gamma}$ is proper (and $p$ is a
homomorphism), this implies that $p|_X$ is proper. Therefore
$X \cap p^{-1}(0)$ is compact. Since $p$ is a homomorphism, we
conclude that $X \cap p^{-1}(0)$ is a compact subgroup of~$X
\cap \real^d$. However,
  \begin{equation} \label{abelian-nocpct}
  \mbox{$\real^d$ has no nontrivial compact subgroups,}
  \end{equation}
  so we conclude that $X \cap p^{-1}(0)$ is trivial, as desired.

\begin{STEP}
  Completion of the proof.
  \end{STEP}
  From Steps~\ref{superabelPf-ponto}
and~\ref{superabelPf-kerp}, and the fact that $X$ is a closed
subgroup of $\real^k \times \real^d$, we see that $X$ is the
graph of a well-defined continuous homomorphism $\hat\varphi
\colon \real^k \to \real^d$.
  Also, because
  $\graph(\varphi) \subset \graph( \hat\varphi)$,
  we know that $\hat\varphi$ extends~$\varphi$.
  \qed \end{proof}

To generalize this proof to the situation where $\integer^k$,
$\real^k$, and~$\real^d$ are replaced by more interesting
solvable groups $\Gamma$, $G$, and~$H$, we need a closed
subgroup~$X$ to substitute for the span of $\hat\Gamma$.
Looking at the proof, we see that the crucial properties
of~$X$ are that it is a connected subgroup that
contains~$\hat\Gamma$ (so $p(X)$ is a connected subgroup
of~$\real^k$ that contains $\operatorname{dom} \varphi$
(see Step~\ref{superabelPf-ponto})), and that
$X/\hat\Gamma$ is compact (see~(\ref{abelian-syndetic})). These
properties are captured in the following definition.

\begin{definition} \label{SyndeticDefn}
  A \emph{syndetic hull}\index{hull!syndetic} of a
subgroup~$\Gamma$ of a Lie group~$G$ is a subgroup~$X$ of~$G$,
such that
  $X$ is connected,
  $X$ contains~$\Gamma$, and
  $X/\Gamma$ is compact.
  \end{definition}

Thus, the same proof applies in any situation where the
following three properties hold:
  \begin{enumerate} \renewcommand{\theenumi}{\alph{enumi}}
  \item \label{needfull}
  no connected, proper subgroup of~$G$ contains~$\Gamma$
(see~(\ref{abelian-fullsubgroup}));
  \item \label{needsyndetic}
  every closed subgroup $\hat\Gamma$ of $G \times H$ has
a syndetic hull (see~(\ref{abelian-syndetic}));
  and
  \item \label{neednocpct}
  $H$ has no nontrivial compact subgroups
(see~(\ref{abelian-nocpct})).
  \end{enumerate}
  Two of these properties pose little difficulty:
  \begin{enumerate}
  \item[({\protect\ref{needfull}})] If $\Gamma$ is a lattice in
a $1$-connected, solvable Lie group~$G$, then no connected,
proper subgroup of~$G$ contains~$\Gamma$
\fullsee{solvable}{full}.
  \item[({\protect\ref{neednocpct}})] If $H$ is a $1$-connected,
solvable Lie group, then $H$ has no nontrivial compact
subgroups \fullsee{solvable}{nocpct}.
  \end{enumerate}
  However, Property~\pref{needsyndetic} may fail, as is
illustrated by the following example.

\begin{example} \label{nosynd}
  Let
  $$ \Gamma =
  \begin{pmatrix}
  1&\integer& \integer\cr
  0&1&0\cr
  0&0&1\cr
  \end{pmatrix}
  \subset G_3
  \mbox{\qquad and\qquad}
  S =
  \bigset{
  \begin{pmatrix}
  1&t& \real\cr
  0&1&0\cr
  0&0&e^{2\pi i t}\cr
  \end{pmatrix}
  }{
  t \in \real
  }
  \; . $$
(see~(\ref{G3Defn})).
  Then $S$ is the only reasonable candidate to be a syndetic
hull of~$\Gamma$ in~$G_3$.
  However, $S$ is not closed under multiplication, so it is not
a subgroup of~$G_3$. Thus, one sees that $\Gamma$ does not have
a syndetic hull in~$G_3$.
  \end{example}

The upshot is that proving superrigidity (in the setting of
solvable groups) reduces to the problem of showing that
syndetic hulls exist. It turns out that Zariski dense
subgroups always have a syndetic hull. (The reader can easily
verify that the subgroup~$\Gamma$ of Example~\ref{nosynd} is
not Zariski dense in~$G_3$.) However, the following key result
(which will be proved in Sect.~\ref{SyndeticSect}) shows that
a much weaker hypothesis suffices: $\closure{\Ad_G \Gamma}$
need only contain a maximal compact subgroup, not all of
$\closure{\Ad_G G}$.

\begin{theorem} \label{syndintro}
  If $\Gamma$ is a closed subgroup of a  connected, solvable
Lie group~$G$, such that
  \begin{equation} \label{MainAssump}
  \mbox{$\closure{\Ad_G \Gamma}$ contains a maximal compact
subgroup of $\closure{\Ad_G G}^\circ$,}
  \end{equation}
  then $\Gamma$ has a syndetic hull\index{hull!syndetic} in~$G$.
Furthermore, if $G$ is $1$-connected, then the syndetic hull is
unique.
  \end{theorem}

For example, if $G$ is a $1$-connected, $\real$-split solvable
group (that is, if $G$ is an upper triangular subgroup of
$\GL(n,\real)$), then $\closure{\Ad_G G}$ has no
compact subgroups, so the hypothesis of the theorem is
trivially satisfied.

\begin{corollary}
  In a $1$-connected, $\real$-split solvable group, syndetic
hulls exist and are unique.
  \end{corollary}

So the proof applies:

\begin{proposition} \label{Saitosuper}
  Suppose
  \begin{itemize}
  \item $G_1$ and $G_2$ are $1$-connected, $\real$-split
solvable groups,
  \item $\Gamma_1$ is a lattice in~$G_1$, and
  \item $\varphi \colon \Gamma_1 \to G_2$ is a homomorphism.
  \end{itemize}
  Then $\varphi$ extends uniquely to a continuous
homomorphism $\hat\varphi \colon G_1 \to G_2$.
  \end{proposition}

\begin{corollary}[Saito 1957] \label{Saitoiso}
  Suppose $\Gamma_i$ is a lattice in a $1$-connected,
$\real$-split solvable group~$G_i$, for $i = 1,2$.
  If $\Gamma_1 \iso \Gamma_2$, then $G_1 \iso G_2$. \end{corollary}

Let us now use Theorem~\ref{syndintro} to prove
a superrigidity theorem.

\begin{proof}[of Theorem~\ref{solvsuper}]
   We are given a homomorphism $\varphi \colon \Gamma \to
\GL(d,\complex)$.

\setcounter{CASE}{0}

\begin{CASE} \label{solvsuperPf-latt}
  Assume $\Gamma$ is a lattice in~$G$.
  \end{CASE}
  Let
  \begin{itemize}
  \item $H = \closure{\varphi(\Gamma)}$,
  \item $\hat G = G \times H$,
  \item
  $\hat\Gamma = \graph(\varphi)
  \subset G \times H $, and
  \item $p \colon G \times H \to G$ be the natural projection
onto the first factor.
  \end{itemize}

We use the proof of Proposition~\ref{ZktoRd}, so there are
only two issues to address. First, we need to show that
$\hat\Gamma$ has a syndetic hull~$X$ in $\hat G$. Second,
because $H$ may not be $1$-connected, we do not have
property~\pref{neednocpct}, the analogue
of~\pref{abelian-nocpct}.

Recall that Zariski closures are virtually connected. (This is
stated formally in Lemma~\ref{Zarvirtconn} below.) Hence, $H$ has
only finitely many components, so, by passing to a finite-index
subgroup of~$\Gamma$, we may assume that $\varphi(\Gamma) \subset
H^\circ$, so $\hat \Gamma \subset G \times H^\circ = \hat
G^\circ$.

By assumption, $\closure{\Gamma}^\circ$ contains a maximal compact
subgroup~$S$ of $\closure{G}^\circ$, and, by definition,
$\closure{\varphi(\Gamma)}^\circ$ contains a maximal compact
subgroup~$T$ of~$H^\circ$. Therefore, the
projection of $\closure{\hat\Gamma}^\circ$ to each factor of
$\closure{G}^\circ \times H^\circ$ contains a maximal compact
subgroup of that factor. However, $\closure{\hat\Gamma}$ is
diagonally embedded in $\closure{G} \times H$, so it probably
does not contain the product $S \times T$, which is a maximal
compact subgroup of $\closure{G}^\circ \times H^\circ$. Thus,
Theorem~\ref{syndintro} probably does not apply directly.
However, $S \times T$ is contained in $\closure{\hat\Gamma}
T$, so the rather technical Theorem~\ref{virtsynd} below,
which can be proved in almost exactly the same way as
Theorem~\ref{syndintro}, does apply. So we conclude that some
finite-index subgroup of $\hat\Gamma$ has a syndetic hull~$X$
in~$\hat G^\circ$, as desired. (Note that, because
$\graph(\hat\varphi) = X$ contains a finite-index subgroup
of~$\hat\Gamma$, the homomorphism~$\hat\varphi$ virtually
extends~$\varphi$.)

Theorem~\ref{virtsynd} asserts that we may take the syndetic
hull~$X$ to be simply connected; thus, $X$ has no nontrivial
compact subgroups. Hence, the subgroup $X \cap p^{-1}(e)$ also
has no compact subgroups. Assumption~(\ref{abelian-nocpct}) was
used only to obtain this conclusion, so we have no need for
\pref{neednocpct}.

\begin{CASE}
  The general case.
  \end{CASE}
  From Theorem~\ref{syndintro}, we know that $\Gamma$ has a
syndetic hull~$B$. So $\Gamma$ is a lattice in~$B$, and, by
assumption, $\closure{\Gamma} = \closure{G} \supset
\closure{B}$. Therefore, Case~\ref{solvsuperPf-latt} implies
that $\varphi$ virtually extends to a continuous homomorphism
$\varphi^* \colon B \to \GL(d,\complex)$.

Now, because $B$ is connected, and $\closure{B} =
\closure{G}$, one can show that $[G,G] \subset B$. So it is
not hard to extend~$\varphi^*$ to a continuous homomorphism
$\hat\varphi \colon G \to \GL(d,\complex)$.
  \qed \end{proof}

\begin{theorem} \label{virtsynd}
  Let $\Gamma$ be a discrete subgroup of a connected, solvable,
linear Lie group~$G$.

If there is a compact subgroup~$S$ of~$\closure{\Gamma}$ and a
compact subgroup~$T$ of~$G$, such that
    $ST$ is a maximal compact subgroup of~$\closure{G}^\circ$,
  then some finite-index subgroup~$\Gamma'$ of~$\Gamma$ has a
simply connected syndetic hull in~$G$.
  \end{theorem}

\begin{literature}
  Definition~\ref{SyndeticDefn} is slightly modified from
\cite{FriedGoldman}. (In our terminology, they proved that
every solvable subgroup~$\Gamma$ of $\GL(\ell,\complex)$
virtually has a syndetic hull in~$\closure{\Gamma}$.)
  Theorem~\ref{syndintro} appears in \cite{Witte-SolvRig}.
  For the special case where $G_1$ and~$G_2$ are nilpotent,
Corollary~\ref{Saitoiso} was proved by Malcev, and this
special case appears in \cite[Theorem~2.11, p.~33]{Raghunathan?}.
  Theorem~\ref{virtsynd} is from \cite{Witte-Sarith}.
  \end{literature}

\section{Solvable Lie Groups and Zariski Closed Subgroups}
  \label{SolvZarSect}

We now recall (without proof) some rather standard
results on solvable Lie groups and Zariski closures.

\subsection{Solvable Lie Groups and Their Subgroups}

\begin{remark} \label{SolvDefnsSame}
  Although the definition of ``solvable'' given in
Defn.~\ref{solv<>triang} is not the usual one, the
Lie--Kolchin Theorem\index{theorem!Lie--Kolchin} \cite[Theorem~17.6,
pp.~113--114]{Humphreys-Algic} implies that a connected
subgroup~$G$ of $\GL(\ell,\complex)$ satisfies
\pref{solv<>triang} if and only if it is solvable\index{group!solvable
Lie}\index{Lie!group!solvable} in the usual sense. Thus, this naive description is
adequate for our purposes.

Also, Ado's Theorem\index{theorem!Ado's} \cite[Theorem~3.18.16,
pp.~246--247]{Varadarajan} implies that every $1$-connected,
solvable Lie group is isomorphic to a closed subgroup of some
$\GL(\ell,\complex)$, so there is no loss of generality in
restricting our attention to linear groups.
  \end{remark}

The following observation is immediate from the usual definition
of solvability:

\begin{lemma} \label{[GG]<G}
  If $G$ is a nontrivial, connected, solvable Lie group, then
  $$\dim [G,G] < \dim G \; .$$
  \end{lemma}

\begin{proposition} \label{solvable}
  Let $H$ be a connected subgroup of a $1$-connected, solvable
Lie group~$G$.
  \begin{enumerate}
  \item \label{solvable-H=Rn}
  $H$ is closed, simply connected, and diffeomorphic to
some~$\real^d$;
  \item \label{solvable-full}
  If $G/H$ is compact, then $H = G$;
  \item \label{solvable-nocpct}
  If $C$ is a compact subgroup of~$G$, then $C$ is trivial.
  \end{enumerate}
  \end{proposition}

\begin{lemma} \label{G/Hsc}
  Let $Q$ be a closed subgroup of a connected, solvable
group~$G$.
  \begin{enumerate}
  \item \label{G/Hsc-->}
  If $G/Q$ is simply connected, then $Q$ is connected, and $Q$
contains a maximal compact subgroup of~$G$;
  \item \label{G/Hsc-<-}
  If $Q$ has only finitely many components, and $Q$ contains a
maximal compact subgroup of~$G$, then $Q$ is connected,
and $G/Q$ is simply connected.
  \end{enumerate}
  \end{lemma}

\begin{lemma} \label{cpctconjsubgrp}
  If $G$ is any Lie group with only finitely many connected
components, then
  \begin{enumerate}
  \item \label{cpctconjsubgrp-exist}
  $G$ has a maximal compact subgroup, and
  \item \label{cpctconjsubgrp-conj}
  all maximal compact subgroups of~$G$ are conjugate to each
other.
  \end{enumerate}
  \end{lemma}

\subsection{Zariski Closed Subgroups of $\GL(\ell,\complex)$}

The following definition formalizes the idea that a subgroup
is Zariski closed if it is defined by polynomial functions.
Also, we are thinking of $\GL(\ell,\complex)$ as being a real
variety of dimension $2\ell^2$, rather than a complex variety
of dimension~$\ell^2$.

\begin{definition}
  A subset $X$ of~$\real^N$ is \emph{Zariski closed}\index{Zariski!closed} 
if there is a (finite or infinite) collection $\{P_k\}$ of
real polynomials in $N$ variables, such that
  $$ X = \{\, x \in \real^N \mid
  \mbox{$P_k(x) = 0$, for all~$k$}\,\}
  \; .$$

Let $\varphi \colon \GL(\ell,\complex) \to \real^{\ell^2+2}$ be
the identification of $\GL(\ell,\complex)$ with a (Zariski
closed) subset of $\real^{\ell^2+2}$ given by listing the real
and imaginary parts of the determinant and of each matrix
entry:
  $$ \varphi(g) =
  \bigl(
  \Re(\det g), \Im(\det g),
  \Re g_{1,1}, \Im g_{1,1},
  \Re g_{1,2}, \Im g_{1,2},
  \ldots,
  \Re g_{\ell,\ell}, \Im g_{\ell,\ell}
  \bigr)
  \; .$$
  A subgroup~$Q$ of $\GL(\ell,\complex)$ is \emph{Zariski
closed}\index{Zariski!closed} if $\varphi(Q)$ is a Zariski closed
subset of $\real^{\ell^2+2}$.

The \emph{Zariski closure}\index{Zariski!closure}~$\closure{\Gamma}$ 
of a subgroup~$\Gamma$ of
$\GL(\ell,\complex)$ is the (unique) smallest Zariski closed
subgroup that contains~$\Gamma$.
  \end{definition}

\begin{lemma} \label{Zarvirtconn}
  Any Zariski closed subgroup has only finitely many connected
components.
  \end{lemma}

\begin{lemma} \label{ZarCN}
  Let $H$ be a Zariski-closed subgroup of $\GL(\ell,\complex)$.
  \begin{enumerate}
  \item \label{ZarCN-C}
  For any subgroup~$\Gamma$ of~$H$, the centralizer
$C_H(\Gamma)$ is Zariski closed;
  \item \label{ZarCN-N}
  For any connected subgroup~$U$ of~$H$, the normalizer
$N_H(U)$ is Zariski closed.
  \end{enumerate}
  \end{lemma}

\begin{corollary} \label{Zarnorm}
  If $G$ is a connected subgroup of
$\GL(\ell,\complex)$, then $\closure{G}$ normalizes~$G$.
  \end{corollary}

\begin{theorem}[Borel Density Theorem] \label{BDT}
  If $\Gamma$ is a closed subgroup of a connected, solvable Lie
group~$G$, such that $G/\Gamma$ is compact, then
  $\closure{\Ad_G \Gamma} T = \closure{\Ad_G G}$, for every
maximal compact subgroup~$T$ of $\closure{\Ad_G
G}^\circ$.\index{theorem!Borel density}\index{density!Borel --- theorem}
  \end{theorem}

\begin{corollary} \label{BDT+T}
  Suppose $\Gamma$ is a closed subgroup of a connected, solvable
Lie group~$G$, such that $G/\Gamma$ is compact. If
\pref{MainAssump} holds, then $\closure{\Ad_G \Gamma} =
\closure{\Ad_G G}$.
  \end{corollary}

\begin{literature}
  Propositions~\fullref{solvable}{H=Rn} and
\fullref{solvable}{nocpct} and Lemma~\ref{cpctconjsubgrp} appear
in \cite[Theorems~12.2.2, 12.2.3, and 15.3.1]{Hochschild-Lie}.
  Proposition~\fullref{solvable}{full} is due to G.~D.~Mostow
\cite[Proposition~11.2]{MostowFSS}.
  Lemma~\ref{G/Hsc} follows from the homotopy long exact sequence
of the fibration $Q \to G \to G/Q$; cf.\
\cite[Lemma~2.17]{Witte-SolvRig}.
  Lemma~\ref{Zarvirtconn} appears in
\cite[Theorem~3.6]{PlatonovRapinchuk}.
  Lemma~\fullref{ZarCN}{C} follows from
\cite[Prop~8.2b]{Humphreys-Algic}.
  Lemma~\fullref{ZarCN}{N} follows from the proof of
\cite[Theorem~3.2.5]{ZimmerBook}.
  A generalization of Theorem~\ref{BDT} appears in
\cite[Corollary~4.2]{Dani-BDT}.
   \end{literature}

\section{Existence of Syndetic Hulls} \label{SyndeticSect}

Constructing a syndetic hull requires some way to show that a
subgroup is connected. The following result on intersections
of connected subgroups is our main tool in this regard.

\begin{proposition} \label{GcapQconn}
  Let $G$ and~$Q$ be solvable Lie subgroups of
$\GL(\ell,\complex)$. If
  \begin{itemize}
  \item $G$ is connected,
  \item $Q$ is Zariski closed {\upshape(}or, more generally,
$Q$ has finite index in~$\closure{Q}${\upshape)}, and
  \item $Q$ contains a maximal compact subgroup of
$\closure{G}^\circ$,
  \end{itemize}
  then $G \cap Q$ is connected.
  \end{proposition}

\begin{proof}
  Let $T$ be a maximal compact subgroup of $\closure{G}^\circ$
that is contained in~$Q$.

\setcounter{CASE}{0}

  \begin{CASE} \label{GcapQconnPf-QinG}
  Assume $Q \subset \closure{G}^\circ$.
  \end{CASE}
  Because $\closure{G}$ normalizes~$G$ \mysee{Zarnorm}, we know
that $Q$ normalizes~$G$, so $GQ$ is a subgroup of
$\GL(\ell,\complex)$. Since $Q$ contains the
maximal compact subgroup~$T$ of~$GQ$, we see that $GQ/Q$ is
simply connected \fullsee{G/Hsc}{<-}. Hence $G/(G \cap Q)
\diffeo GQ/Q$ is simply connected, so $G \cap Q$ is connected
\fullsee{G/Hsc}{->}.

\begin{CASE}
  The general case.
  \end{CASE}
  Because $\closure{G}^\circ \cap Q$ contains the maximal
compact subgroup~$T$ of $\closure{G}^\circ$,
Case~\ref{GcapQconnPf-QinG} implies that
  $$G \cap Q
  = (G \cap \closure{G}^\circ) \cap Q
  = G \cap (\closure{G}^\circ \cap Q)
  $$
  is connected, as desired.
  \qed \end{proof}

The following corollary is obtained by using the proposition
to show that $\hat{G} \cap \hat{Q}$ is connected, where
$\hat{G} = \graph(\rho)$ and $\hat{Q} = \closure{G} \times Q$.

\begin{corollary} \label{invimgConn}
  Let
  \begin{itemize}
  \item $G$ be a connected, solvable Lie group,
  \item $\rho\colon G\to \GL(d,\complex)$ be a
finite-dimensional, continuous representation, and
  \item $Q$ be a Zariski closed subgroup
of~$\GL(d,\complex)$, such that
  $Q$ contains a maximal compact subgroup
of~$\closure{\rho(G)}^\circ$.
  \end{itemize}
  Then $\rho^{-1}(Q)$ is connected.
  \end{corollary}

\begin{corollary} \label{allconn}
  Let $\Gamma$ be a closed subgroup of a $1$-connected, solvable
Lie group~$G$, such that \pref{MainAssump} holds.
  \begin{enumerate}
  \item \label{allconn-Z}
  If $\Gamma \subset Z(G)$, then $Z(G)$ is connected;
  \item \label{allconn-C}
  If $\Gamma$ is abelian, then $C_G(\Gamma)$ is connected;
  \item \label{allconn-N}
  If $U$ is  any connected subgroup of~$G$ that is normalized
by~$\Gamma$, then $N_G(U)$ is connected.
  \end{enumerate}
  \end{corollary}

\begin{proof}
  To simplify the notation, let us assume $\closure{\Gamma}$
contains a maximal compact subgroup of~$\closure{G}^\circ$,
ignoring the adjoint representation. (Without this
simplification, the proof would use
Corollary~\ref{invimgConn}, with $\rho = \Ad_G$, instead of
using Proposition~\ref{GcapQconn}, as we do here.) Under this
assumption, Proposition~\ref{GcapQconn} implies that if $Q$ is
any Zariski closed subgroup of $\GL(\ell,\complex)$ that
contains~$\Gamma$, then $G \cap Q$ is connected. We apply this
fact with:

\pref{allconn-Z} $Q = C_{\closure{G}}(G)$ \fullsee{ZarCN}{C};
so $Z(G) = C_G(G) = G \cap Q$ is connected.

\pref{allconn-C} $Q = C_{\closure{G}}(\Gamma)$
\fullsee{ZarCN}{C}; so $C_G(\Gamma) = G \cap Q$ is connected.

\pref{allconn-N} $Q = N_{\closure{G}}(U)$ \fullsee{ZarCN}{N};
so $N_G(U) = G \cap Q$ is connected.
  \qed \end{proof}

\begin{theorem} \label{syndexist}
  Suppose $\Gamma$ is a closed subgroup of a $1$-connected,
solvable Lie group~$G$. If \pref{MainAssump} holds, then
$\Gamma$ has a unique syndetic hull in~$G$.
  \end{theorem}

\begin{proof}
  Let us first prove uniqueness: suppose $S_1$ and~$S_2$ are
syndetic hulls of~$\Gamma$. We have $\closure{\Ad_G \Gamma} =
\closure{\Ad_G S_i}$ for $i = 1,2$ \mysee{BDT+T}; so
$\closure{\Ad_G S_1} = \closure{\Ad_G S_2}$. Therefore $S_1$
and~$S_2$ normalize each other \mysee{Zarnorm}, so $S_1 S_2$ is
a subgroup of~$G$. It is simply connected
\fullsee{solvable}{H=Rn}, and $S_1 S_2/S_2 \diffeo S_1/(S_1
\cap S_2)$ is compact (because $\Gamma \subset S_1 \cap S_2$),
so Lemma~\fullref{solvable}{full} implies that $S_2 = S_1 S_2$;
thus $S_1 \subset S_2$. Similarly, $S_2 \subset S_1$.
Therefore $S_1 = S_2$, so the syndetic hull, if it exists, is
unique.

We now prove existence. We may assume, by induction on $\dim
G$ \mysee{[GG]<G}, that $\Gamma \cap [G,G]$ has a unique syndetic
hull~$U$ in $[G,G]$ (note that $\closure{\Ad_G [G,G]}$ is
unipotent, so it has no nontrivial compact subgroups).

\setcounter{CASE}{0}

\begin{CASE} \label{syndexistPf-Gabel}
  Assume $G$ is abelian.
  \end{CASE}
  Because $\Gamma$ is a normal subgroup of~$G$, we may consider
the quotient group $G/\Gamma$: let $K/\Gamma$ be a maximal
compact subgroup of~$G/\Gamma$. By definition, $K/\Gamma$ is
compact. Also, $G/K \diffeo (G/\Gamma)/(K/\Gamma)$ is simply
connected \fullsee{G/Hsc}{<-}, so $K$ is connected
\fullsee{G/Hsc}{->}. Therefore $K$ is a syndetic hull
of~$\Gamma$.

\begin{CASE}  \label{syndexistPf-center}
  Assume $\Gamma \subset Z(G)$.
  \end{CASE}
  Because $Z(G)$ is connected \fullsee{allconn}{Z}, we know,
from Case~\ref{syndexistPf-Gabel}, that $\Gamma$ has a
syndetic hull~$S$ in $Z(G)$. Then $S$ is a syndetic hull
of~$\Gamma$ in~$G$.

\begin{CASE} \label{syndexistPf-abelGamma}
  Assume $\Gamma$ is abelian.
  \end{CASE}
  We have $\Gamma \subset C_G(\Gamma)$, and $C_G(\Gamma)$ is
connected \fullsee{allconn}{C}, so, from
Case~\ref{syndexistPf-center}, we know that $\Gamma$ has a
syndetic hull~$S$ in $C_G(\Gamma)$. Then $S$ is also a
syndetic hull of~$\Gamma$ in~$G$.

\begin{CASE} \label{syndexistPf-Unormal}
  Assume $U$ is a normal subgroup of~$G$.
  \end{CASE}
  Let $\tau \colon G \to G/U$ be the natural homomorphism.
Then, because $U/(\Gamma \cap U)$ is compact, we see that
$\Gamma U$ is closed, so $\tau(\Gamma)$ is a closed subgroup of
$\tau(G)$. Also, because
$[\Gamma,\Gamma] \subset U$, we have
  $$ [\tau(\Gamma), \tau(\Gamma)]
  = \tau \bigl( [\Gamma,\Gamma] \bigr)
  \subset \tau(U)
  = e
  \; ,$$
  so $\tau(\Gamma)$ is abelian. Thus, from
Case~\ref{syndexistPf-abelGamma}, we know that $\tau(\Gamma)$
has a syndetic hull~$S$ in $\tau(G)$. Then
$\tau^{-1}(S)$ is a syndetic hull of~$\Gamma$ in~$G$.

\begin{CASE}
  The general case.
  \end{CASE}
  The uniqueness of the syndetic hull~$U$ implies that $\Gamma$
normalizes~$U$; that is, $\Gamma \subset N_G(U)$. Now $N_G(U)$
is connected \fullsee{allconn}{N}, so, from
Case~\ref{syndexistPf-Unormal}, we know that $\Gamma$ has a
syndetic hull~$S$ in $N_G(U)$; then $S$ is also a syndetic
hull of~$\Gamma$ in~$G$.
  \qed \end{proof}

\begin{corollary} \label{syndexistnotsc}
  Suppose $\Gamma$ is a closed subgroup of a connected, solvable
Lie group~$G$. If \pref{MainAssump} holds, then $\Gamma$ has a
syndetic hull in~$G$.
  \end{corollary}

\begin{proof}
  Write $G = \cover G/Z$ and $\Gamma = \cover \Gamma/Z$, where
$Z$ is some discrete, normal subgroup of the center of the
universal cover~$\cover G$ of~$G$. If $S$ is any syndetic hull
of~$\cover \Gamma$, then $S/Z$ is a syndetic hull of~$\Gamma$.
  \qed \end{proof}

\begin{remark} If $G$ is not simply connected, then
  syndetic hulls may not be unique.
   (For example, $e$ and $\torus$ are two syndetic hulls of~$e$
in~$\torus$.)
  \end{remark}

Proposition~\ref{GcapQconn} has the following corollary.
Theorem~\ref{virtsynd} is proved almost exactly the same way
as Theorem~\ref{syndexist}, but using this corollary in place
of Proposition~\ref{GcapQconn}. However, a small additional
argument is needed when $G$ is abelian, to show that the
syndetic hull can be chosen to be simply connected in this
base case.

\begin{corollary}
\label{GcapQvirtconn}
  Let $G$ and~$Q$ be solvable Lie subgroups of
$\GL(\ell,\complex)$. If
  \begin{itemize}
  \item $G$ is connected,
  \item $Q$ is Zariski closed, and
  \item there are compact subgroups~$S$ of~$Q$ and
$T$~of~$\closure{G}$, such that $ST$ is a maximal compact
subgroup of $\closure{G}^\circ$,
  \end{itemize}
  then $G \cap Q$ is virtually connected.
  \end{corollary}

\begin{proof}
  By replacing $Q$ with $\closure{G} \cap Q$, we may assume $Q
\subset \closure{G}$. Also, by replacing $T$ with a subgroup,
we may assume $S \cap T$ is finite.
  From the structure theory of solvable Zariski closed
subgroups \cite[Theorems~19.3 and 34.3b]{Humphreys-Algic}, we have
  $\closure{G}^\circ = (ST) \ltimes V$
  and $Q^\circ = S \ltimes (Q \cap V)$,
  where $V$ is the
subgroup generated by the elements of~$\closure{G}$, all of
whose eigenvalues are real and positive;
  then, because $Q$ contains the maximal compact subgroup~$S$
of $SV = \bigl( \closure{G \cap SV} \bigr)^\circ$,
Proposition~\ref{GcapQconn} implies
  $Q \cap \bigl( G \cap (SV) \bigr)^\circ$ is connected. This
is a finite-index subgroup of $Q \cap G$, because $Q^\circ
\subset SV$, and $G \cap (SV)$ is virtually connected.
Therefore, $Q \cap G$ is virtually connected.
  \qed \end{proof}

\vskip1cm
\noindent
{\it Acknowledgments.}   An anonymous referee provided helpful comments.

\bibliographystyle{alpha}

\clearpage
\flushbottom
\printindex

\end{document}